\documentclass[12pt]{amsart}
\usepackage{amssymb,amsmath}
\usepackage{geometry}    
\usepackage{mathrsfs}

\usepackage{graphicx}

\usepackage{vmargin}
\usepackage{tikz}                   

\newtheorem{theorem}{Theorem}[section]

\begin{document}

\title[Construction of special $3$-systems]{On a question
of Schmidt and Summerer concerning $3$-systems}

\author{Johannes Schleischitz} 

\thanks{Middle East Technical University, Northern Cyprus Campus, Kalkanli,
	G\"uzelyurt \\
	johannes.schleischitz@univie.ac.at}


\begin{abstract}
Following a suggestion of W.M. Schmidt and L. Summerer,	
we construct a proper $3$-system $(P_{1},P_{2},P_{3})$
with the property $\overline{\varphi}_{3}=1$. 
In fact, our method
generalizes to provide $n$-systems with $\overline{\varphi}_{n}=1$, for
arbitrary $n\geq 3$. We visualize our constructions with graphics.
We further present explicit examples of numbers $\xi_{1}, \ldots, \xi_{n-1}$ that
induce the $n$-systems in question.
\end{abstract}

\maketitle

{\footnotesize{

Math Subject Classification 2010: 11J13, 11H06}}

\vspace{1mm}

\section{Parametric Diophantine approximation
in dimension three} \label{intro}

Let $\xi_{1},\xi_{2}$ be real numbers so that the set
$\{ 1,\xi_{1},\xi_{2}\}$ is linearly independent over $\mathbb{Q}$.
For $q>0$ a parameter,
let $K(q)$ be the box of points 
$(z_{0},z_{1},z_{2})\in\mathbb{R}^{3}$ that satisfy
\[
\vert z_{0}\vert \leq e^{2q},\qquad \vert z_{1}\vert \leq e^{-q},
\qquad \vert z_{2}\vert \leq e^{-q}.
\] 
Further let $\Lambda$ be the lattice consisting of the points
$\{ (x,\xi_{1}x-y_{1}, \xi_{2}x-y_{2}): x,y_{1},y_{2}\in\mathbb{Z}\}$. 
The successive minima $\lambda_{1}(q),\lambda_{2}(q),\lambda_{3}(q)$
of $K(q)$ with respect to
$\Lambda$ as functions of $q$ contain the essential
information on the simultaneous rational approximation 
to $\xi_{1},\xi_{2}$. 
It is convenient to study the logarithms of the
functions $\lambda_{j}(q)$, denoted by
$L_{j}(q)=\log \lambda_{j}(q)$ for $j=1,2,3$.
These functions have the nice property that their slopes are
among $\{-2,1\}$, and their sum is absolutely bounded
uniformly in the parameter $q$.
These properties motivated
Schmidt and Summerer~\cite{s2} to define so called $3$-systems.
A $3$-system $P=(P_{1},P_{2},P_{3})$ is a triple of
functions $P_{j}: [0,\infty)\to\mathbb{R}$ 
with slopes among $\{-2,1\}$ with the properties that
$P_{1}(0)=P_{2}(0)=P_{3}(0)=0$, 
$P_{1}(q)\leq P_{2}(q)\leq P_{3}(q)$ and
$P_{1}(q)+P_{2}(q)+P_{3}(q)=0$ for every $q\geq 0$. Hence,
locally in a neighborhood of any $q>0$, 
precisely one of the three functions decays
while the other two rise, unless 
$q$ is a switch point where some 
$P_{j}$ are not differentiable (change slope). 
Moreover, for $P$ to be a $3$-system,
it is additionally required that if at a switch point
$q$ some $P_{i}$ changes 
from falling to rising and some other $P_{j}$ from rising to 
falling, then $i<j$ unless $P_{i}(q)=P_{j}(q)$. 
It has been shown in~\cite{s2}
that every function triple $(L_{1},L_{2},L_{3})$ as above,
associated to some $(\xi_{1},\xi_{2})$,
corresponds to a $3$-system $P$ up to a bounded amount,
and conversely by Roy~\cite{roy} that for any $3$-system
$P$ there exist $\xi_{1},\xi_{2}$ satisfying the 
$\mathbb{Q}$-linear independence condition above
and so that $\sup_{q>0}\max_{j=1,2,3}\vert P_{j}(q)-L_{j}(q)\vert \ll 1$. Roy's result employs a minor technical condition on 
the mesh of the system $P$, we do not rephrase it here.
Both results~\cite{s2,roy} are established in more generality.

For given $\xi_{1},\xi_{2}$ with induced funtions $L_{j}(q)$,
let $\varphi_{j}(q)=L_{j}(q)/q$ and put 
 \[
 \underline{\varphi}_{j}= \liminf_{q\to \infty} \varphi_{j}(q),
 \qquad\qquad  \overline{\varphi}_{j}= \limsup_{q\to \infty} \varphi_{j}(q),
 \]
 for $j=1,2,3$. Since $L_{j}$ have slopes $-2$ and $1$ only, it is clear that
\begin{equation} \label{eq:clear}
-2\leq \underline{\varphi}_{j}\leq \overline{\varphi}_{j}\leq 1,
\qquad\qquad j=1,2,3.
\end{equation}
 By
 virtue of the results from~\cite{roy, s2} quoted above, 
 in the sequel we will identify the
 values $\underline{\varphi}_{j}, \overline{\varphi}_{j}$
 with quantities derived from an associated $3$-system
 $P$ via 
 \begin{equation} \label{eq:id}
 \underline{\varphi}_{j}\longleftrightarrow \liminf_{q\to\infty} \frac{P_{j}(q)}{q},
 \qquad
  \overline{\varphi}_{j}\longleftrightarrow  \limsup_{q\to\infty} \frac{P_{j}(q)}{q},
  \qquad\qquad j=1,2,3, 
 \end{equation}
and vice versa. 
M. Laurent~\cite{laurent} provided estimates
for classical exponents of approximation related to
any pair $(\xi_{1},\xi_{2})$ that is $\mathbb{Q}$-linearly
independent with $\{ 1\}$. As pointed out in~\cite{s4}
they translate
into the language of the functions $\varphi_{j}$ as  
 \begin{align}
 0\leq \underline{\varphi}_{3}\leq 
 \overline{\varphi}_{3}&\leq 1, \nonumber \\ 
 \underline{\varphi}_{3}+\underline{\varphi}_{3}\overline{\varphi}_{1}+\overline{\varphi}_{1}&=0, \label{eq:0} \\
 2\underline{\varphi}_{1}+\overline{\varphi}_{3}&\leq
 -\underline{\varphi}_{3}(3+\underline{\varphi}_{1}+2\overline{\varphi}_{3}), \label{eq:1} \\
  2\overline{\varphi}_{3}+\underline{\varphi}_{1}&\geq
 -\overline{\varphi}_{1}(3+\overline{\varphi}_{3}+2\underline{\varphi}_{1}).  \label{eq:2}
 \end{align}
Schmidt and Summerer~\cite{s4} recently provided
additional information by including the second successive minimum
in the picture.

\begin{theorem}[Schmidt/Summerer, 2017] \label{ss}
For any $\xi_{1},\xi_{2}$ with
$\{ 1,\xi_{1},\xi_{2}\}$ linearly independent over $\mathbb{Q}$,
if $0\leq \underline{\varphi}_{3}<1$,
additionally to the above relations we have
\begin{equation} \label{eq:eq1}
\overline{\varphi}_{2}\leq \overline{\Omega}:=
\frac{\overline{\varphi}_{1}-\underline{\varphi}_{1}}{2-\overline{\varphi}_{1}-\overline{\varphi}_{1}\underline{\varphi}_{1}},
\end{equation}
and 
\begin{equation} \label{eq:eq2}
\underline{\varphi}_{2}\geq \underline{\Omega}:=
\frac{\underline{\varphi}_{3}-\overline{\varphi}_{3}}{2-\underline{\varphi}_{3}-\overline{\varphi}_{3}\underline{\varphi}_{3}}.
\end{equation} 
Moreover, these estimates are best possible in the sense that
for given numbers $\underline{\varphi}_{1},\overline{\varphi}_{1},\underline{\varphi}_{3},\overline{\varphi}_{3}$ with
$0\leq \underline{\varphi}_{3}<1$ and \eqref{eq:0}, \eqref{eq:1}, \eqref{eq:2}
there are $\xi_{1},\xi_{2}$ 
with $\{ 1,\xi_{1},\xi_{2}\}$ linearly independent over $\mathbb{Q}$
for whose approximation constants we have
$\underline{\varphi}_{2}=\underline{\Omega}$
and $\overline{\varphi}_{2}=\overline{\Omega}$.
\end{theorem}

Schmidt and Summerer enclose a remark to Theorem~\ref{ss}
pointing out that in the case $\underline{\varphi}_{3}= 1$ excluded
in its claim, we have $\overline{\varphi}_{2}=1$ and
$\overline{\varphi}_{1}=\underline{\varphi}_{2}= -1/2$ (by
mistake they denoted $-1/3$ instead in~\cite{s4}).   
However, there is a gap in
Theorem~\ref{ss} concerning the {\em existence} of graphs with the property
$\underline{\varphi}_{3}=1$, and related real numbers
$\xi_{1}, \xi_{2}$. In~\cite{s4} they state
''But one really should prove that
$\underline{\xi}=(\xi_{1},\xi_{2})$ with 
$(1, \xi_{1},\xi_{2})$ linearly independent over 
$\mathbb{Q}$ with $\overline{\varphi}_{3}= 1$ exist. We
invite the reader to construct a proper 3-system $P$ with this property.''
The main purpose of this paper is to provide the desired construction. 
Before we turn to constructing the $3$-system, we point out
that explicit examples of $\mathbb{Q}$-linearly independent
 $\{1,\xi_{1},\xi_{2}\}$ inducing
$\underline{\varphi}_{3}=\overline{\varphi}_{3}= 1$ can
be derived from previous results of the author. Concretely~\cite[Corollary~2.11]{j1}, 
upon putting $k=n-1=2$ and $C=\infty$, yields the following example.

\begin{theorem} \label{hres}
Let
\begin{equation} \label{eq:given}
\xi_{1}= \sum_{k=1}^{\infty} 10^{-(2k-1)!}, \qquad\qquad 
\xi_{2}= \sum_{k=1}^{\infty} 10^{-(2k)!}.
\end{equation}
Then
\begin{align}
\underline{\varphi}_{1}&=-2, \qquad \underline{\varphi}_{2}=-\frac{1}{2}, \qquad 
\underline{\varphi}_{3}=1, \label{eq:h3} \\
\overline{\varphi}_{1}&=-\frac{1}{2}, \qquad \overline{\varphi}_{2}=1, \qquad\quad 
\overline{\varphi}_{3}=1.  \label{eq:h4} 
\end{align} 
\end{theorem}

While the results in~\cite{j1}
are originally formulated in the language of another type of
exponents,
the two types of exponents determine each other
via the identities of~\cite[Theorem~1.4]{s1},
and we derive Theorem~\ref{hres}. 
We note that for the sole purpose
of $\overline{\varphi}_{3}=1$, as desired
in~\cite{s4} and rephrased above, in fact any
numbers $\xi_{1}, \xi_{2}$ which are simultaneously approximable
to any order by rational numbers can be chosen. 
In particular, one may choose
the pair $(\xi,\xi^{2})$ with $\xi$ any Liouville number, 
see~\cite[Theorem 3.1]{glasnik}. However, then
we always have $\underline{\varphi}_{3}=0$.
For Liouville's constant given as $\xi=10^{-1!}+10^{-2!}+10^{-3!}+\cdots$,
by~\cite[Theorem~3.2]{glasnik} in place
of \eqref{eq:h3}, \eqref{eq:h4} we have
\begin{align}
\underline{\varphi}_{1}&=-2, \qquad \underline{\varphi}_{2}=-\frac{1}{2}, \qquad 
\underline{\varphi}_{3}=0 \label{eq:h1} \\
\overline{\varphi}_{1}&=0, \qquad\quad \overline{\varphi}_{2}=1, \qquad\quad 
\overline{\varphi}_{3}=1.  \label{eq:h2} 
\end{align}

Alternatively to the above examples,
the pure existence of pairs $(\xi_{1},\xi_{2})$ inducing
$\overline{\varphi}_{3}=1$ (or $\underline{\varphi}_{3}=1$)
also follows from Roy's results~\cite{roy} and
\cite[Theorem~11.5]{roy2} (the latter result, already quoted in~\cite{s4},
provides an explicit description of the spectrum of sixtuples $\underline{\varphi}_{1},\underline{\varphi}_{2},\underline{\varphi}_{3},\overline{\varphi}_{1},\overline{\varphi}_{2},\overline{\varphi}_{3}$
by a system of complicated inequalities). The main concern of the
question of Schmidt and Summerer appears to be the construction
of a suitable $3$-system, carried out in 
Section~\ref{c} below.

\section{Construction of a $3$-system with $\overline{\varphi}_{3}=1$}

We want to 
present an effective construction
of a $3$-system with \eqref{eq:h3}, \eqref{eq:h4},
in particular $\overline{\varphi}_{3}=1$.
It resembles
the combined graph $(L_{1},L_{2},L_{3})$ with respect to 
the pair $(\xi_{1}, \xi_{2})$ 
in \eqref{eq:given}, in an idealized form.
In fact the resulting $3$-system 
can be interpreted
as the idealized extremal case of the regular graph defined 
in~\cite{s3}, for the parameter $\rho=\infty$.
In Section~\ref{c2} we will briefly sketch how to modify the method
to obtain a graph with \eqref{eq:h1}, \eqref{eq:h2} instead, and give
generalizations to $n$-systems.

\subsection{The construction} \label{c}

We construct the graphs piecewise as follows.
Let 
\[
0<l_{0}< l_{1}< l_{2}< l_{3}<\cdots,
\]
be a fast increasing lacunary sequence
of real numbers with the property 
\begin{equation} \label{eq:assuan}
\lim_{i\to\infty} \frac{l_{i+1}}{l_{i}}=\infty.
\end{equation}
Let $r_{0}=0$. In the interval $[r_{0},l_{0}]=[0,l_{0}]$ 
let $P_{1}$ decay
with slope $-2$ and $P_{2},P_{3}$ rise with slope $1$,
so that $P_{1}(l_{0})=-2l_{0}$ and 
$P_{2}(l_{0})=P_{3}(l_{0})=l_{0}$. Let $w_{0}=l_{0}$
for consistency with later notation.
Let $l_{0}$ be the first switch point where $P_{1}$ starts to rise
and $P_{2}$ starts to decay. Then the graph of $P_{1}$ will
meet the graph of $P_{2}$ at some point $(r_{1},P_{1}(r_{1}))$
with $r_{1}>l_{0}$.
We may assume $l_{1}>r_{1}$.
In the interval $[r_{1}, l_{1}]$ we define $P_{1}$ as decaying
with slope $-2$ again and the other two functions rising with
slope $1$. Note that $P_{3}(l_{1})=l_{1}$ since it has not 
changed slope yet. Assume this construction of the graphs in 
$[0,l_{1}]$
was step 0 of our construction. Now we carry out
how to complete the process with identical
steps 1,2,3,... where in step $i$ we define the graphs 
of $P_{1},P_{2},P_{3}$ in the
interval $[l_{i},l_{i+1}]$.
At position $q=l_{1}$ we let $P_{1}$ and
$P_{3}$ change slopes so that $P_{1}$ rises with slope $1$
and $P_{3}$ decays
with slope $-2$. The function $P_{2}$ still rises with slope $1$. We keep these slopes
until $P_{2}$ meets $P_{3}$ at position $q=w_{1}$. 
Then we let $P_{2}$ decay with
slope $-2$ and the other functions rise with slope $1$ until
$P_{2}$ meets $P_{1}$ at some point $(r_{2},P_{1}(r_{2}))$. 
We may assume $l_{2}>r_{2}$.
Then we let $P_{1}$ decay
with slope $-2$ up to $q=l_{2}$, and the other
two functions rise with slope $1$ in this interval.
This completes step 1. 
At $q=l_{2}$ we let $P_{1}$ again switch
from decaying to rising and conversely for $P_{3}$, and so on.
When we repeat the whole
process ad infinitum, 
we claim that $P_{1},P_{2},P_{3}$ represent 
the combined graph
of a $3$-system with the 
properties \eqref{eq:h3}, \eqref{eq:h4}. A sketch of such a $3$-system
in an initial interval is shown in Figure 1 below. For size reasons we used
the slopes $-1,1/2$ instead of $-2,1$, thereby sketching $P_{j}(q)/2$ for $j=1,2,3$.

\begin{tikzpicture}

\draw[thick,->] (0,0) -- (13,0) node[anchor=west] {q};
\draw[thick,->] (0,-6) -- (0,6) node[anchor=south] {P(q)};

\draw (0,0) -- (0.5,0.25);
\draw (0,0) -- (0.5,-0.5) node[anchor=north west] {$P_{1}$};

\draw (0.5,-0.5) -- (3,0.75);
\draw (0.5,0.25) -- (3,-2.25) node[anchor=north west] {$P_{1}$};
\draw (0.5,0.25) -- (3,1.5) node[anchor=south west] {$P_{3}$};
\draw (3,1.5) -- (10,-5.5) node[anchor=north west] {$P_{1}$};
\draw (3,0.75) -- (10,4.25) node[anchor=south west] {$P_{3}$};
\draw (3,-2.25) -- (10,1.25);
\draw (10,4.25) -- (12.5,1.75) node[anchor=north west] {$P_{2}$};
\draw (10,1.25) -- (12.5,2.5) node[anchor=south west] {$P_{3}$};
\draw (10,-5.5) -- (12.5,-4.25) node[anchor=west] {$P_{1}$};

\fill[gray] (0.5,0) circle (0.06cm) node[anchor=north] {$l_{0}$};
\fill[gray] (3,0) circle (0.06cm) node[anchor=north] {$l_{1}$};
\fill[gray] (10,0) circle (0.06cm) node[anchor=north] {$l_{2}$};
\fill[gray] (5.5,0) circle (0.06cm) node[anchor=north] {$r_{2}$};
\fill[gray] (3.5,0) circle (0.06cm) node[anchor=north] {$w_{1}$};
\fill[gray] (12,0) circle (0.06cm) node[anchor=north] {$w_{2}$};
\fill[gray] (1,0) circle (0.06cm) node[anchor=south] {$r_{1}$};

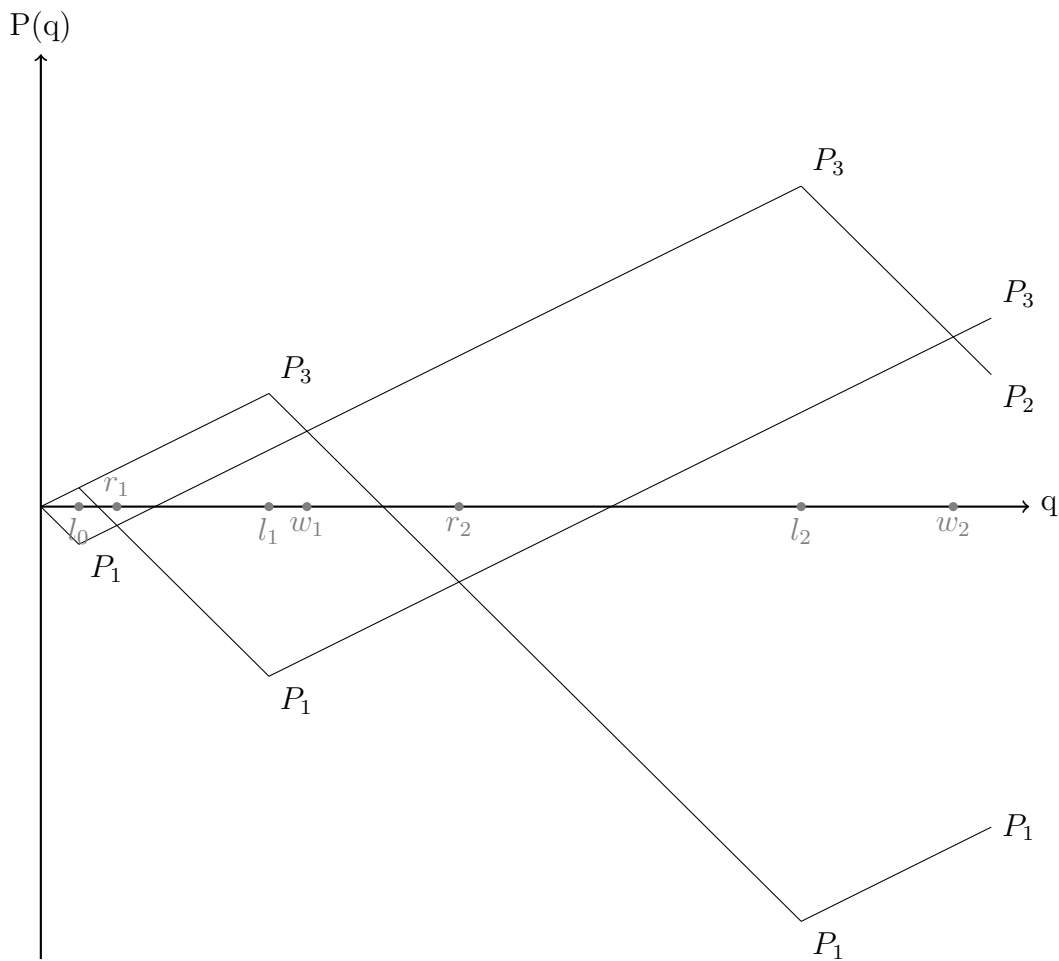
\node at (5.5,-6.5) {Figure 1: Visualization of case $\overline{\varphi}_{3}=1$,
slopes scaled by factor $1/2$};

\end{tikzpicture}

\subsection{The proof}
Keep in mind for the following that
the switching positions in our construction are ordered
\[
0=r_{0}<l_{0}=w_{0}<r_{1}<l_{1}<w_{1}<r_{2}<l_{2}<w_{2}<\cdots,
\]
and also the identification \eqref{eq:id}.
First it is clear that the process yields the combined 
graph of a $3$-system $P$. Indeed, by construction there is always
precisely one $P_{j}$ decaying, there
are infinitely many positions where $P_{1}=P_{2}$ 
and $P_{2}=P_{3}$ respectively hold, and the switches 
occur in a way that respects the additional $3$-system condition
on a local maximum having higher index than a local
minimum at switch points
mentioned in the introduction. To obtain \eqref{eq:h3}, \eqref{eq:h4},
we first look at positions $q=l_{i}$ and claim that 
\begin{equation} \label{eq:cl}
\lim_{i\to\infty} \frac{P_{1}(l_{i})}{l_{i}}=-2,
\qquad\qquad \lim_{i\to\infty} \frac{P_{2}(l_{i})}{l_{i}}= \lim_{i\to\infty}\frac{P_{3}(l_{i})}{l_{i}}=1.
\end{equation}
By the identification \eqref{eq:id} and by \eqref{eq:clear}
this implies $\underline{\varphi}_{1}=-2$ and 
$\overline{\varphi}_{2}=\overline{\varphi}_{3}=1$.
By construction $P_{1}$ decays with slope $-2$ in intervals
of the form $I_{t}:=[r_{t},l_{t}]$ for $t\geq 0$ and rises in intervals
$J_{t}:=[l_{t-1},r_{t}]$ for $t\geq 1$. We next check that
\begin{equation} \label{eq:rl}
r_{t}<2l_{t-1}, \qquad \qquad t\geq 1.
\end{equation}
We trivially have $P_{3}(l_{t-1})-P_{1}(l_{t-1})\leq l_{t-1}-(-2l_{t-1})=3l_{t-1}$. 
On the other hand,
since $P_{1}$ decays in $J_{t}$ with slope $-2$
whereas $P_{3}$ rises with slope $1$, the function $P_{3}-P_{1}$
has slope $3$ in $J_{t}$ so that they must meet within
distance $3l_{t-1}/3=l_{t-1}$ in the first coordinate
on the right from $l_{t-1}$. This
intersection point has first coordinate $r_{t}$, and we deduce
\eqref{eq:rl}.   

The estimate \eqref{eq:rl} and the assumption
\eqref{eq:assuan} clearly imply that the sums of the 
lengths of the intervals $I_{t}$ over 
$t=1,2,\ldots,i$ exceeds
the according sums of the intervals $J_{t}$ by any given factor
$\rho>0$ for large enough $i$, i.e.
\[
\sum_{t=0}^{i} \vert I_{t}\vert > \rho \sum_{t=1}^{i} \vert J_{t}\vert, \qquad\qquad i\geq i_{0}(\rho).
\]
Thus since
\[
l_{i}= \sum_{t=0}^{i} \vert I_{t}\vert+\sum_{t=1}^{i} \vert J_{t}\vert
\]
and
\[
P_{1}(l_{i})= -2 \sum_{t=0}^{i} \vert I_{t}\vert+\sum_{t=1}^{i} \vert J_{t}\vert,
\]
indeed for sufficiently large $i$ we have
\[
\frac{P_{1}(l_{i})}{l_{i}} < -\frac{2+\rho^{-1}}{1+\rho^{-1}}.
\]
As we can choose $\rho$ arbitrarily large indeed
 $\lim_{i\to\infty} P_{1}(l_{i})/l_{i}=-2$, hence
$\underline{\varphi}_{1}=-2$ by \eqref{eq:clear}, \eqref{eq:id}. Since $P_{2}$ 
and $P_{3}$ rise
with slope $1$ in any $I_{t}$ we infer the 
remaining claims of \eqref{eq:cl} by a very similar argument,
or directly by using the bounded sum property at $q=l_{i}$.

Next we show
\begin{equation} \label{eq:next}
\lim_{q\to\infty} \frac{P_{3}(q)}{q}=1.
\end{equation}
By construction $P_{3}$ has local minima precisely
at positions $w_{i}$ and it rises with slope $1$ everywhere outside
of the intervals $[l_{i},w_{i}]$, in which it decays with slope $-2$. 
In view of \eqref{eq:clear} it suffices to
check that 
\begin{equation} \label{eq:thus}
\liminf_{i\to\infty} \frac{P_{3}(w_{i})}{w_{i}}\geq 1.
\end{equation}
By construction
\[
P_{3}(w_{i})= l_{0} - 2\sum_{j=0}^{i} (w_{j}-l_{j}) +\sum_{j=0}^{i-1} (l_{j+1}-w_{j}).
\]
Hence, in view of \eqref{eq:assuan},
to verify \eqref{eq:thus} it suffices to check
\begin{equation} \label{eq:holz}
\lim_{i\to\infty} \frac{w_{i}}{l_{i}}=1. 
\end{equation}
Now by construction
in the interval $[l_{i},w_{i}]$ the function $P_{2}$ rises with slope $1$ whereas
$P_{3}$ decays with slope $-2$, hence 
$w_{i}=l_{i}+u_{i}$ with $u_{i}$ defined implicitly by the identity 
$P_{3}(l_{i})-2u_{i}= P_{2}(l_{i})+ u_{i}$, that 
is $w_{i}=l_{i}+(P_{3}(l_{i})-P_{2}(l_{i}))/3$. On the other hand, by
\eqref{eq:cl} we have $P_{2}(l_{i})= l_{i}(1+o(1))$ and $P_{3}(l_{i})= l_{i}(1+o(1))$, hence inserting we derive $w_{i}=l_{i}(1+o(1))$ as $i\to\infty$, as desired. 
Thus \eqref{eq:next} is shown.

Finally we show that 
\begin{equation} \label{eq:cl2}
\lim_{i\to\infty} \frac{P_{1}(r_{i})}{r_{i}}=\lim_{i\to\infty} \frac{P_{2}(r_{i})}{r_{i}}=-\frac{1}{2}.
\end{equation}
Since by construction the local maxima of $P_{1}$ 
and the local minima of $P_{2}$ both are attained precisely at
the positions $r_{i}$, the remaining identities
from \eqref{eq:h3} and \eqref{eq:h4} are implied.
Let $K_{t}=[w_{t-1},r_{t}]$, so that $K_{t}\subseteq J_{t}$ and by 
\eqref{eq:assuan}, \eqref{eq:holz} the complement $J_{t}\setminus K_{t}$ is small
compared to $J_{t}$.
In $K_{t}$, the function $P_{1}$ rises with slope $1$ 
whereas $P_{2}$ decays with slope $-2$. Moreover,
by \eqref{eq:cl} and \eqref{eq:holz} and since the slopes are bounded
\[
\lim_{i\to\infty} \frac{P_{1}(w_{i})}{w_{i}}=-2,
\qquad\qquad \lim_{i\to\infty} \frac{P_{2}(w_{i})}{w_{i}}=1.
\] 
Combining these two facts and by definition of $r_{i}$,
for large $i$ we readily conclude $r_{i}= w_{i-1}(2-o(1))$ and
thus the
asymptotic value at $r_{i}$ is $P_{1}(r_{i})=P_{1}(w_{i-1})+r_{i}-w_{i-1}= w_{i-1}(-1+o(1))$, hence indeed $P_{1}(r_{i})/r_{i}=-1/2+o(1)$ for large $i$.
Thus \eqref{eq:cl2} holds and
the proof is finished.

\subsection{Generalizations and variations} \label{c2}
A similar construction as in Section~\ref{c} can be 
done in arbitrary dimension $n$, where
the slopes of the $P_{j}$ are among
$\{-n+1,1\}$. Instead of one sequence $(w_{i})_{i\geq 0}$
with $l_{i}<w_{i}<r_{i+1}$, we obtain $n-2$ sequences
$(w_{i}^{h})_{i\geq 0}, 1\leq h\leq n-2$, induced by
positions where $P_{n-h+1}$ meets $P_{n-h}$, ordered
$l_{i}<w_{i}^{1}<w_{i}^{2}<\cdots<w_{i}^{n-2}<r_{i+1}$. 
We derive $n$-systems 
$P=(P_{1},\ldots,P_{n})$  
whose approximation constants (via identification \eqref{eq:id}) satisfy
\[
\underline{\varphi}_{1}=-n+1, \qquad \underline{\varphi}_{2}=\frac{2-n}{2}, \qquad
\underline{\varphi}_{j}=1, \quad 3\leq j\leq n,
\]
and
\[
\overline{\varphi}_{1}=\frac{2-n}{2}, \qquad \qquad
\overline{\varphi}_{j}=1, \quad 2\leq j\leq n.
\]
Again this resembles the special case $\rho=\infty$ of the
regular graph~\cite{s3} in dimension $n$, and
suitable numbers $(\xi_{1},\ldots,\xi_{n-1})$ 
inducing these approximation constants arise
from~\cite[Corollary~2.11]{j1} upon taking $k=n-1, C=\infty$, a
particular choice is
\[
\xi_{j}= \sum_{k=0}^{\infty} 10^{-(k(n-1)+j)!}, \qquad\qquad 1\leq j\leq n-1.
\] 

Finally, we sketch the construction of a $3$-system $P$ with
the properties \eqref{eq:h1}, \eqref{eq:h2} in place
of \eqref{eq:h3}, \eqref{eq:h4}. We have
to alternate between the construction of
Section~\ref{c} and another type of intermediate
construction. Take $i$ a large integer and follow the construction from
Section~\ref{c} up to $q=q_{0}=:l_{i}$. Recall 
$P_{1}(l_{i})\approx -2l_{i}$ and 
$P_{j}(l_{i})\approx l_{i}$ for $j=2,3$ by \eqref{eq:cl}. Then we make
the first intermediate construction. 
Starting from $q_{0}$, let $P_{1}$ rise
with slope $1$ and $P_{2},P_{3}$ decay with slope roughly $-1/2$ in not too
short intervals. The latter can be easily realized by changing the slopes
of $P_{2}, P_{3}$ rapidly so that there are many positions $q$ with 
equality $P_{2}(q)=P_{3}(q)$.
One may take these equality positions an arithmetic sequence
$b_{0},b_{1}=b_{0}+D,b_{2}=b_{0}+2D,\ldots,b_{h}=b_{0}+hD$ with
some $b_{0}\geq q_{0}, h\geq 0$ and some small increment $D>0$,
in the following way.
Fix $D>0$ small.
Let $(b_{0},P_{2}(b_{0}))$ be the intersection point of the line
passing through $(q_{0},P_{2}(q_{0}))$ with slope $1$ (graph of $P_{2}$)
and the line passing through $(q_{0},P_{3}(q_{0}))$ with 
slope $-2$ (graph of $P_{3}$), corresponding to $w_{i}$ in Section~\ref{c}. 
In $[b_{0}, b_{0}+D/2]$,
let $P_{2}$ decay with slope $-2$ and $P_{3}$ rise with slope $1$. 
Then at $q_{0}+D/2$ interchange the slopes, such that
at $b_{1}=b_{0}+D$ we have $P_{2}(b_{1})= P_{3}(b_{1})= P_{2}(b_{0})-D/2$.  
We repeat this procedure and stop at the largest index $h$ so that
the resulting graphs of $P_{2},P_{3}$ remain positive on $[0,b_{h}]$.
For simplicity let $\tilde{q}:=b_{h}$.
Notice that
$P_{j}(b_{l})-P_{j}(b_{0})=-(b_{l}-b_{0})/2=-lD/2$ for $l=0,1,\ldots,h$.
Therefore, by \eqref{eq:cl} and since $D$ is small, it is easy to see that 
$\vert P_{j}(\tilde{q})\vert$ are all small for $j=1,2,3$. 
Now starting at $\tilde{q}$,
let $P_{1}, P_{3}$ rise with slope $1$ and $P_{2}$ decay with
slope $-2$ until the graphs of $P_{1}$ and $P_{2}$ meet at some position $q_{1}$. 
Since $\vert P_{j}(\tilde{q})\vert$ are all small, 
the expressions $q_{1}-\tilde{q}$ 
and $\vert P_{j}(q_{1})\vert$ for $j=1,2,3$, are small (like $o(q_{1})$) as well. 
This ends the first intermediate construction, illustrated in Figure 2 below (again slopes are
scaled with factor $1/2$).
Now we essentially apply the initial
construction (step $0$)
from Section~\ref{c} from the interval $[0,l_{1}]$ 
again,
starting from $q=q_{1}$ instead of $q=0$. Let us denote by $q_{2}$
the right endpoint in this construction, that is
the value corresponding to $l_{1}$ from Section~\ref{c}.
Notice that
$P_{1}$ has a local minimum inside the interval $[q_{1}, q_{2}]$,
corresponding to $l_{0}$ from Section~\ref{c}, and another one at
the right endpoint $q_{2}$.
Since $\vert P_{j}(q_{1})\vert$ 
are small for $j=1,2,3$, the $P_{j}$ indeed behave in $[q_{1},q_{2}]$
essentially like they do in the construction
of Section~\ref{c} in the interval $[0,l_{1}]$ (see Figure~1).
In particular, as for $q=q_{0}$, at $q=q_{2}$ again 
we have $P_{1}(q_{2})\approx -2q_{2}$ and 
$P_{j}(q_{2})\approx q_{2}$ for $j=2,3$. Hence at this point we
again switch to the intermediate construction to define the $P_{j}$ 
in some interval $[q_{2},q_{3}]$. We
repeat this iterative process of constructing $P$ in
$[q_{2k},q_{2k+1}]$ and
then in $[q_{2k+1},q_{2k+2}]$, for all $k\geq 1$. 
It can be checked that the resulting combined graph
satisfies \eqref{eq:h1}, \eqref{eq:h2}. Notice hereby that the condition
$\underline{\varphi}_{2}=-1/2$ forced us to copy the behavior of the
$P_{j}$ on $[0,l_{1}]$, and not only on $[0,l_{0}]$, in intervals $[q_{2k+1},q_{2k+2}]$.
The procedure can again be 
generalized to dimension $n$ to provide $n$-systems with the properties
\[
\underline{\varphi}_{1}=-n+1, \qquad \underline{\varphi}_{2}=\frac{2-n}{2}, \qquad
\underline{\varphi}_{j}=0, \quad 3\leq j\leq n,
\]
and
\[
\overline{\varphi}_{1}=0, \qquad\qquad
\overline{\varphi}_{j}=1, \quad 2\leq j\leq n.
\]

\begin{tikzpicture}

\draw[thick,->] (0,0) -- (13.3,0) node[anchor=west] {q};
\draw[thick,->] (0,-6) -- (0,6) node[anchor=south] {P(q)};

\draw (2.7,-5.2) -- (3,-5.5)  node[anchor=north west] {$P_{1}$} ;
\draw (2.7,3.35) -- (3,3.5)   node[anchor=south west] {$P_{3}$};
\draw (3,2) -- (2.7,1.85)    node[anchor=north east] {$P_{2}$};

\draw (3,3.5) -- (4,2.5);
\draw (3,2) -- (4,2.5);
\draw (3,-5.5) -- (13.2,-0.4);

\draw (4,2.5) -- (4.5,2.75);
\draw (4,2.5) -- (4.5,2);

\draw (4.5,2.75) -- (5,2.25);
\draw (4.5,2) -- (5,2.25);

\draw (5, 2.25) -- (5.5, 2.5);
\draw (5, 2.25) -- (5.5,1.75);

\draw (5.5,2.5) -- (6,2);
\draw (5.5,1.75) -- (6,2);

\draw (6,2) -- (6.5,2.25);
\draw (6,2) -- (6.5,1.5);

\draw (6.5,2.25) -- (7,1.75);
\draw (6.5,1.5) -- (7,1.75);

\draw (7,1.75) -- (7.5,2);
\draw (7,1.75) -- (7.5,1.25);

\draw (7.5,2) -- (8,1.5);
\draw (7.5,1.25) -- (8,1.5);

\draw (8,1.5) -- (8.5,1);
\draw (8,1.5) -- (8.5,1.75);

\draw (8.5,1) -- (9,1.25);
\draw (8.5,1.75) -- (9,1.25);

\draw (9,1.25) -- (9.5,1.5);
\draw (9,1.25) -- (9.5,0.75);

\draw  (9.5,1.5) -- (10,1);
\draw  (9.5,0.75) -- (10,1);

\draw (10,1) -- (10.5,1.25);
\draw (10,1) -- (10.5,0.5);

\draw (10.5,1.25) -- (11,0.75);
\draw (10.5,0.5) -- (11,0.75);

\draw (11,0.75) -- (11.5,1);
\draw (11,0.75) -- (11.5,0.25);

\draw (11.5,1) -- (12,0.5);
\draw (11.5,0.25) -- (13.2,1.1);

\draw (12,0.5) -- (13.2,-0.7);

\fill[gray] (13,0) circle (0.06cm) node[anchor=north] {$q_{1}$};
\fill[gray] (12,0) circle (0.06cm) node[anchor=north] {$b_{8}=\tilde{q}$};
\fill[gray] (3,0) circle (0.06cm) node[anchor=north] {$q_{0}=l_{i}$};
\fill[gray] (4,0) circle (0.06cm) node[anchor=north] {$b_{0}$};
\fill[gray] (4,0) circle (0.06cm) node[anchor=north] {$b_{0}$};
\fill[gray] (5,0) circle (0.06cm) node[anchor=north] {$b_{1}$};
\fill[gray] (6,0) circle (0.06cm) node[anchor=north] {$b_{2}$};
\fill[gray] (7,0) circle (0.06cm) node[anchor=north] {$b_{3}$};
\fill[gray] (8,0) circle (0.06cm) node[anchor=north] {$b_{4}$};
\fill[gray] (9,0) circle (0.06cm) node[anchor=north] {$b_{5}$};
\fill[gray] (10,0) circle (0.06cm) node[anchor=north] {$b_{6}$};
\fill[gray] (11,0) circle (0.06cm) node[anchor=north] {$b_{7}$};

\draw[dashed] (4,2.5) -- (12,0.5);

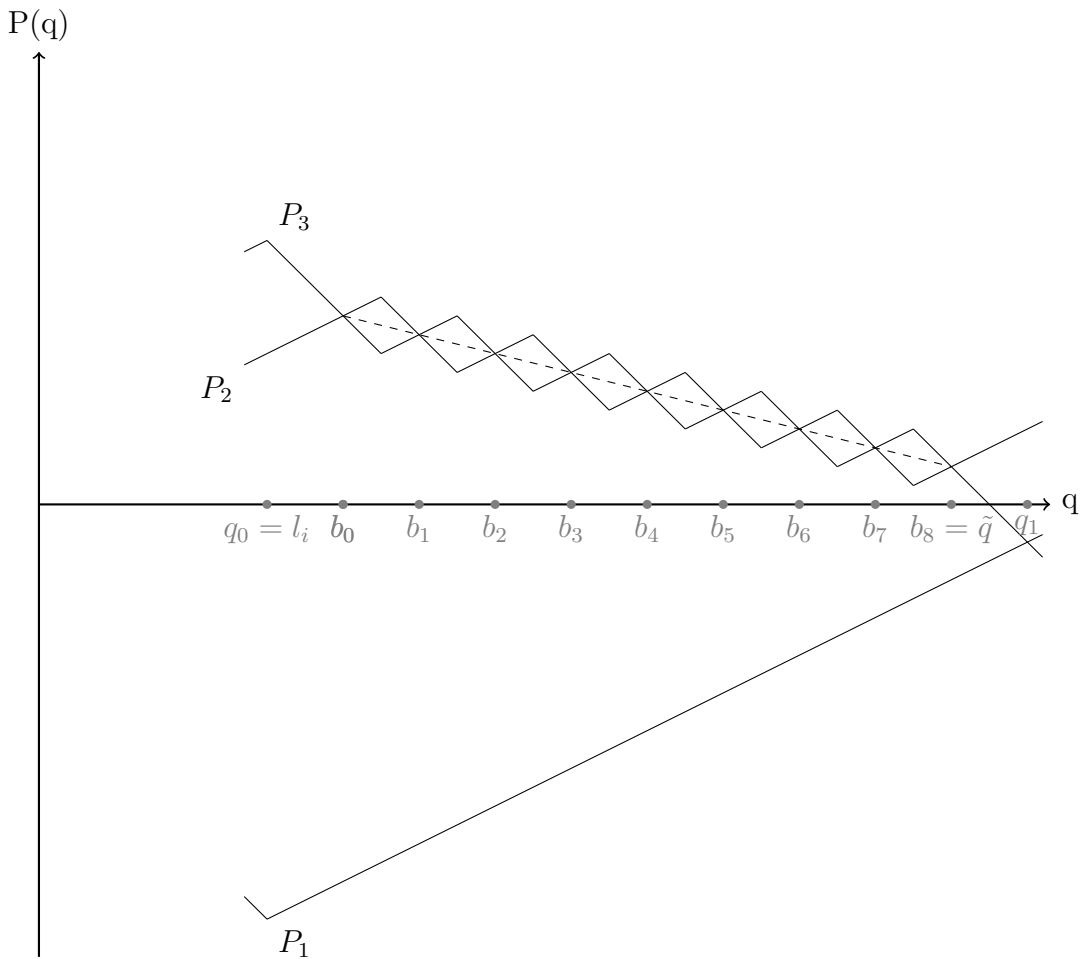
\node at (6.3,-7) {Figure 2: Intermediate construction in $[q_{0}, q_{1}]$,
slopes scaled by factor $1/2$};

\end{tikzpicture}

\vspace{0.5cm}

{\em The author thanks the referee for the careful reading and for
	pointing out some
inaccuracies! }

\end{document}